\input amstex
\input amsppt.sty
\magnification=1200
\def\Span{\operatorname{Span}}
\def\wn{\operatorname{int}}
\def\epsilon{\varepsilon}
\def\id{\operatorname{id}}
\def\phi{\varphi}
\document

\topmatter
\title
On Bergman completeness of pseudoconvex Reinhardt domains
\endtitle
\author
W\l odzimierz Zwonek
\endauthor
\address
Instytut Matematyki, Reymonta 4, 30-059 Krak\'ow, Poland
\endaddress
\address
current addresse: Carl von Ossietzky Universit\"at Oldenburg,
Fachbereich 6 -- Mathematik, Postfach 2503, 26111 Oldenburg,
Germany
\endaddress
\email
zwonek\@im.uj.edu.pl (zwonek\@mathematik.uni-oldenburg.de)
\endemail
\thanks
The author is a fellow of the Alexander von Humboldt Foundation. The author
has been partially supported by the KBN grant No. 2 P033A 017 14.
\endthanks
\abstract
In the paper we give a precise description 
of Bergman complete bounded pseudoconvex Reinhardt domains
\endabstract
\endtopmatter

\subheading{0. Introduction} The study of the boundary behaviour 
of the Bergman kernel has a very comprehensive literature. 
Very closely related to this problem is the
problem of Bergman completeness of a domain. It is well-known that 
any Bergman complete domain is pseudoconvex (see \cite{Bre}). 
The converse is not true (take the punctured unit disc in $\Bbb C$).
The problem which pseudoconvex domains are Bergman complete 
has a long history (see e.g. \cite{Ohs}, \cite{Jar-Pfl~1}). 
The most general
result, stating that any bounded hyperconvex domain is Bergman complete,
has been recently obtained independently by Z. B\l ocki and P. Pflug
(see \cite{B\l o-Pfl}) and by Z. Herbort (see \cite{Her}). 
The result generalizes
earlier results on Bergman completeness. An example of a bounded
non-hyperconvex domain, which is Bergman complete can be found already 
in dimension one (see \cite{Chen}). Another example of such a domain is given in
\cite{Her}; the example is a bounded pseudoconvex Reinhardt 
domain in $\Bbb C^2$. Motivated by the last example we solve in the paper
entirely the problem, which bounded Reinhardt domains are Bergman complete. 
Since the concept of hyperconvexity in the class of bounded
Reinhardt domains is completely understood (see \cite{Car-Ceg-Wik} and
\cite{Zwo~2}) the paper gives us complete answer about the mutual
relationship between hyperconvexity and the Bergman completeness 
in the considered class of domains. The condition which satisfy
the Bergman complete bounded pseudoconvex Reinhardt domains is expressed
in simple geometrical properties of some convex cones 
associated to Reinhardt pseudoconvex domains. 

The paper is the next one in understanding the completeness with respect to
invariant distances in the class 
of Reinhardt pseudoconvex domains (see \cite{Pfl~2}, \cite{Fu},
\cite{Zwo~1} and \cite{Zwo~2}). 

\subheading{1. Definitions and statement of results}
Let us denote by $E$ the unit disc in $\Bbb C$. By $A_+$ (repsectively,
$A_-$) we
denote the non-negative (respectively, non-positive) numbers from $A$. Put
also $A_*^n:=(A_*)^n$, where $A_*:=A\setminus\{0\}$.

Let $D$ be a 
domain in $\Bbb C^n$. Let us denote by
$L_h^2(D)$ square integrable holomorphic functions on $D$. $L_h^2(D)$
is a Hilbert space. Let
$\{\phi_j\}_{j\in J}$ ($\#J\leq\aleph_0$) be an orthonormal basis of
$L_h^2(D)$.
Then we define
$$
K_D(z):=\sum_{j\in J}|\phi_j(z)|^2,\;z\in D.
$$
Let us call $K_D$ {\it the Bergman kernel of $D$}.
For domains $D$ such that for any $z\in D$ there is $f\in L^2_h(D)$
with $f(z)\neq 0$ (for example for $D$ bounded) we have
$$
K_D(z)=\sup\{\frac{|f(z)|^2}{||f||_{L^2(D)}^2}:f\in L^2_h(D),f\not\equiv 0\}
$$
If $D$ is such that $K_D(z)>0$, $z\in D$ then
$\log K_D$ is a smooth plurisubharmonic function. 
In this case we define
$$
\beta_D(z;X):=\left(\sum_{j,k=1}^n\frac{\partial^2\log K_D(z)}
{\partial z_j\partial \bar z_k}X_j\bar X_k\right)^{1/2},\;
z\in D,\;X\in\Bbb C^n.
$$ 
$\beta_D$ is a pseudometric called {\it the Bergman pseudometric}.

For $w,z\in D$ we put
$$
b_D(w,z):=
\inf\{L_{\beta_D}(\alpha)\}
$$
where the infimum is taken over piecewise $C^1$-curves
$\alpha:[0,1]\mapsto D$ joining $w$ and $z$ and 
$L_{\beta_D}(\alpha):=\int_0^1\beta_D(\alpha(t);\alpha^{\prime}(t))dt$.

We call $b_D$ {\it the Bergman pseudodistance of $D$}.

The Bergman distance (as well as the Bergman metric) is
invariant with respect to biholomorphic mappings.
In other words, for any biholomorphic mapping $F:D\mapsto G$ 
($D,G\subset\subset \Bbb C^n$) we have
$$
b_G(F(w),F(z))=b_D(w,z),\quad\beta_G(F(w);F^{\prime}(w)X)=\beta_D(w;X),
\;w,z\in D,X\in\Bbb C^n.
$$
A bounded domain $D$ is called
{\it Bergman complete } if any $b_D$-Cauchy sequence is convergent to
some point in $D$ with respect to the standard topology of $D$.

As already mentioned any bounded Bergman complete domain is pseudoconvex 
(see \cite{Bre}). 
Let us recall that a bounded domain $D$ is called {\it hyperconvex }
if it admits a continuous negative plurisubharmonic exhaustion function.
Now we may formulate the following very general result:

\proclaim{Theorem 1 {\rm (see \cite{B\l o-Pfl}, \cite{Her})}} 
Let $D$ be a bounded hyperconvex domain in $\Bbb C^n$.
Then $D$ is Bergman complete.
\endproclaim
It is known that the converse implication in Theorem 1 does not hold
(see \cite{Chen} and \cite{Her}). Our aim is to find a precise description
of Bergman complete bounded Reinhardt pseudoconvex domains. Let us underline
that we make no use of Theorem 1.
Before we can formulate the result let us recall 
some standard definitions and results on pseudoconvex Reinhardt domains 
and let us introduce the notions necessary for formulating our theorem.

\vskip2ex

A domain $D\subset\Bbb C^n$ is called {\it Reinhardt } if
$(e^{i\theta_1}z_1,\ldots,e^{i\theta_n}z_n)\in D$ for any
$z\in D$, $\theta_j\in\Bbb R$, $j=1,\ldots,n$.

For a point $z\in\Bbb C_*^n$ we put $\log|z|:=(\log|z_1|,\ldots,\log|z_n|)$.
We denote $\log D:=\{\log|z|:z\in D\cap\Bbb C_*^n\}$.

Let us denote:
$$
\gather
V_j:=\{z\in\Bbb C^n:z_j=0\},\;j=1,\ldots,n;\\
V_I:=V_{j_1}\cap\ldots\cap V_{j_k},\text{ where } I=\{j_1,\ldots,j_k\},\;
1\leq j_1<\ldots<j_k\leq n.
\endgather
$$
We have the following description of pseudoconvex Reinhardt domains.

\proclaim{Proposition 2 {\rm (see \cite{Vla}, 
\cite{Jak-Jar})}} Let $D$ be a Reinhardt domain. Then
$D$ is a pseudoconvex Reinhardt domain if and only if
$$
\gather
\log D \text{ is convex and
for any $j\in\{1,\ldots,n\}$}\\
\text{if $D\cap V_j\neq\emptyset$  and ($z^{\prime},z_j,
z^{\prime\prime})\in D$ then
$(z^{\prime},\lambda z_j,z^{\prime\prime})\in D$
for any $\lambda\in \bar E$}.
\endgather
$$
\endproclaim
From the above description we have following properties.
Assume that $D$ is a Reinhardt pseudoconvex domain and assume that
for some $j\in\{1,\ldots,n\}$ we have $D\cap V_j\neq\emptyset$. Then for
the mapping
$$
\pi_j:D\owns z\mapsto (z_1,\ldots,z_{j-1},0,z_{j+1},\ldots,z_n)\in V_j
$$
we have the property
$$
\pi_j(D)=D\cap V_j.
$$
In particular, $\pi_j(D)$ is a Reinhardt pseudoconvex domain in $\Bbb C^{n-1}$
(after trivial identification). We may go further and we may
formulate the following result.

Assume that $D\cap V_{I}\neq\emptyset$,
$I=\{j_1,\ldots,j_k\}$, $1\leq j_1<\ldots<j_k\leq n$, $k<n$.
Define the mapping
$(\pi_I(z))_j:=0$ if $j\in I$ and $z_j$ otherwise. Then
$\pi_I(D)=D\cap V_{I}$ and $\pi_I(D)$
is a Reinhardt pseudoconvex domain in $\Bbb C^{n-k}$.

\vskip2ex

In view of the above considerations
convex domains in $\Bbb R^n$ play an important role
in the study of pseudoconvex Reinhardt domains.
It turns out that while considering different classes of holomorphic
functions the special role is played by cones associated
to the logarithmic image of the domain. 

We say that $C\subset\Bbb R^n$ is a {\it cone  with vertex at $a$ } if
for any $v\in C$ we have $a+tv\in C$ whenever $t>0$. If, in the sequel,
we do not specify the vertex of the cone, then we shall mean a cone with
vertex at $0$.

Following \cite{Zwo~3} for a convex domain $\Omega\subset \Bbb R^n$ 
and a point $a\in\Omega$
let us  define
$$
\frak C(\Omega,a):=\{v\in\Bbb R^n:a+\Bbb R_+v\subset\Omega\}.
$$
It is easy to see that
$\frak C(\Omega,a)$ is a closed convex cone (with vertex at $0$).
Note that
$$
\frak C(\Omega,a)=\bigcup_{C+a\subset\Omega, C\text{ -- cone}}C=
\text{the largest cone contained in $(\Omega-a)$}.
$$
Moreover,
$\frak C(\Omega,a)=\frak C(\Omega,b)$ for any $a,b\in\Omega$. Therefore,
we may well define $\frak C(\Omega):=\frak C(\Omega,a)$ for some (any)
$a\in\Omega$. 

Note that assuming $0\in \Omega$
we have that $\frak C(\Omega)=h^{-1}(0)$, where $h$ is the Minkowski
functional of $\Omega$.
It is also easy to see that
$$
\frak C(\Omega)=\{0\} \text{ if and only if } \Omega\subset\subset \Bbb R^n.
$$

\vskip2ex

For a pseudoconvex Reinhardt domain $D\subset \Bbb C^n$ we define
$\frak C(D):=\frak C(\log D)$. Let us also define
$$
\align
\tilde\frak C(D)&:=\{v\in\frak C(D):
\overline{\exp(a+\Bbb R_+v)}\subset D\},\\
\frak C^{\prime}(D)&:=\frak C(D)\setminus \tilde\frak C(D),
\endalign
$$
where $a$ is some point from $\log D$.

Let us remark that the definition of $\tilde\frak C(D)$
(and, consequently, that of $\frak C^{\prime}(D)$)
does not depend on the choice of $a\in\log D$ (exactly as in the case of
$\frak C(D)$). It follows easily from the above mentioned
properties of pseudoconvex Reinhardt domains.

Our aim is the following result:
\proclaim{Theorem 3} For a bounded pseudoconvex Reinhardt domain
$D$ in $\Bbb C^n$ the following two conditions are equivalent:
$$
\gather
D \text{ is Bergman complete},\tag{i}\\
\frak C^{\prime}(D)\cap\Bbb Q^n=\emptyset.\tag{ii}
\endgather
$$
\endproclaim 
The existence of vectors from $\frak C^{\prime}(D)$, which are rational
(equivalently are from $\Bbb Z^n$) means that we may embed a
punctured disc in $D$ in such a way that the mapping extends to a mapping
defined on the whole disc and the extension maps $0$ to a point
from the boundary (actually, this mapping is a monomial mapping whose
powers come from the element of $\frak C^{\prime}(D)\cap\Bbb Z^n$). 

In our study the key role will be played by the following criterion:

\proclaim{Theorem 4 {\rm (see \cite{Kob})}} 
Let $D$ be a bounded domain.
Assume that there is a subspace $\Cal E\subset L_h^2(D)$ with
$\bar \Cal E=L_h^2(D)$ such that for any $f\in\Cal E$,
$z^0\in\partial D$ and for any sequence $\{z^{\nu}\}_{\nu=1}^{\infty}
\subset D$
$z^{\nu}\to z^0$ there is a subsequence $\{z^{\nu_j}\}$ such that
$$
\frac{|f(z^{\nu_j})|}{\sqrt{K_D(z^{\nu_j})}}\to 0\text{ as $j\to\infty$}.
\tag{KC}
$$
Then $D$ is Bergman complete.
\endproclaim

As already mentioned 
the first non-hyperconvex pseudoconvex Reinhardt domain, which
is Bergman complete was given by Z. Herbort (see \cite{Her}). His example
was the one with $\frak C(D)=\tilde\frak C(D)=\{0\}\times\Bbb R_-$
and his proof was based on the Kobayashi Criterion applied to some
subspace defined with the help of weight functions. 

The case $n=2$ in Theorem 3 was
proven in \cite{Zwo~3}. In the present paper we use the same methods 
(however, simplified and generalized)
as those given there;
especially we verify \thetag{KC} for a linear subspace
of finite linear combinations of monomials square-integrable on $D$.
It is the simplest possible
subspace, which may be tested in the condition \thetag{KC}. 

For $\alpha=(\alpha_1,\ldots,\alpha_n)\in\Bbb Z^n$ put 
$z^{\alpha}:=z_1^{\alpha_1}\cdot\ldots\cdot z_n^{\alpha_n}$ for
$z\in\Bbb C^n$ such that $z_j\neq 0$ if $\alpha_j<0$. Then we define
$$
\Cal E:=\Span\{z^{\alpha}: z^{\alpha}\in L_h^2(D)\}.
$$
We know that $\bar\Cal E=L_h^2(D)$. 

Note that in order to verify the property
\thetag{KC} at some $z^0$ for $\Cal E$ it is sufficient to show
that this property holds for all functions $z^{\alpha}\in L_h^2(D)$.

\subheading{2. Auxiliary results} 
Below we give some results concerning the description of 
monomial mappings in pseudconvex Reinhardt domains and some results
on diophantine approximation in convex cones, which will be crucial in our
considerations.
Lemmas 5--7 come from \cite{Zwo~3}; nevertheless,
for the sake of completeness we give their proofs here.
\proclaim{Lemma 5 {\rm(cf. \cite{Zwo~3})}} 
Let $D$ be a pseudoconvex Reinhardt domain.
Let $\alpha\in \Bbb Z^n$. Then
$$
z^\alpha\in L_h^2(D) \text{ if and only if }
\langle \alpha+\bold 1,v\rangle<0
\text{ for any } v\in\frak C (D),\; v\neq 0.
$$
\endproclaim
\demo{Proof} Assume that $a=(1,\ldots,1)\in D$.
First we prove the following

{\bf Claim.} Assume that $\frak C(D)\neq\{0\}$.
Then for any $\epsilon>0$ there is a cone $T$
such that $(\log D)\setminus T$ is bounded and

if $v\in T$, $||v||=1$ then there exists $w\in \frak C(D)$ such that
$||w||=1$ and $||v-w||<\epsilon$.

\demo{Proof of the Claim} Let $h$ be a Minkowski functional of $\log D$.
$\log D$ is convex, so
$h$ is continuous. Recall that $h^{-1}(0)=\frak C(D)$. 
From the continuity of $h$
we get that for any $\epsilon>0$ there is $\delta>0$ such that
$\{w\in\Bbb R^n:h(w)\leq\delta,||w||=1\}\subset
\{w\in\Bbb R^n:||w||=1$
and there is $v\in \frak C(D)$, $||v||=1$, $||w-v||<\epsilon\}$.

Now take the cone $T$ to be the smallest cone
containing the set
$\{w\in\Bbb R^n:h(w)\leq\delta,||w||=1\}$. 
Note that $(\log D)\setminus T$ is bounded.
If it were not the case, then there would be $x_{\nu}\to\infty$
such that $x_{\nu}\in(\log D)\setminus T$,
so $h(x_{\nu})<1$, consequently $h(\frac{x_{\nu}}{||x_{\nu}||})<
\frac{1}{||x_{\nu}||}$, so $x_{\nu}\in T$ for $\nu$ large enough
-- contradiction.
\qed
\enddemo

If $\frak C(D)=\{0\}$ then the result is trivial.
Assume that $\frak C(D)\neq\{0\}$. Fix $\alpha\in \Bbb Z^n$ such that 
$z^{\alpha}\in L^2_h(D)$.
Let $v\in \frak C(D)$, $v\neq 0$. We may assume that $|v_n|=1$.
There is an open bounded set $U\subset \Bbb R^{n-1}$ such that
$0\in U\times\{0\}$ and $U\times\{0\}+\Bbb R_+v\subset\log D$.
We have
$$
\multline
\infty>\int_D|z^{\alpha}|^{p}=
\int_{D\cap\Bbb C_*^n}|z^{\alpha}|^{2}=
(2\pi)^n\int_{\log D}e^{2\langle\alpha+\bold 1,x\rangle}
dx_1\ldots dx_n\geq\\
(2\pi)^n\int_0^{\infty}
\big(\int_{U\times\{0\}+x_nv}e^{2\langle\alpha+\bold 1,
x\rangle}dx_1\ldots dx_{n-1}\big)dx_n=
M\int_0^{\infty}e^{2x_n\langle\alpha+\bold 1,v\rangle}dx_n,
\endmultline
$$
from which we get the desired inequality $\langle\alpha+\bold 1,v
\rangle<0$.

Assume now that $\langle\alpha+\bold 1,v\rangle<0$
for any $v\in \frak C(D)$, $v\neq 0$.
Then there is some $\delta>0$ such that
$\langle\alpha+\bold 1,v\rangle\leq-\delta$ 
for any $v\in \frak C(D)$,
$||v||=1$.
Now using Claim we get the existence of a cone $T$ fulfilling
among others the following inequality:
$$
\langle\alpha+\bold 1,v\rangle\leq -\delta/2,\; v\in T,\;||v||=1
$$
or
$$
\langle\alpha+\bold 1,v\rangle\leq(-\delta/2)||v||,\;v\in T.
$$
It follows from the description of $T$ ($(\log D)\setminus T$ is bounded)
that 
$$
\int_{\log D}e^{2\langle\alpha+\bold 1,x\rangle}dx<\infty
\text{ if and only if }\int_{T}
e^{2\langle\alpha+\bold 1,x\rangle}dx<\infty.
$$
And now let us estimate the last expression
$$
\int_{T}
e^{2\langle\alpha+\bold 1,x\rangle}dx\leq
\int_Te^{-\delta||x||}dx\leq
\int_{\Bbb R^n}
e^{-\delta||x||}dx<\infty,
$$
which finishes the proof of the lemma.
\qed
\enddemo
\proclaim{Lemma 6 {\rm(cf. \cite{Zwo~3})}}
Let $H$ be a $k$--dimensional vector subspace of $\Bbb R^n$
such that $H\cap \Bbb Q^n=\{0\}$. Let $\{v^1,\ldots,v^k\}$
be a vector base of $H$. Then the set
$$
\{(\langle\alpha,v^1\rangle,\ldots,\langle\alpha,v^k\rangle):
\alpha\in\Bbb Z^n\}
$$
is dense in $\Bbb R^k$.
\endproclaim
\demo{Proof} Certainly, $k<n$. 
It is easy to see that there is a vector subspace
$\tilde H\supset H$ of dimension $n-1$
such that $\tilde H\cap \Bbb Q^n=\{0\}$. Therefore,
we lose no generality assuming that $k=n-1$.

Moreover, we lose no generality assuming that for a matrix
$$
\tilde V:=
\bmatrix
v^1_1& \ldots& v^{n-1}_1\\
\cdot& \ldots& \cdot \\
v^1_{n-1}&\ldots& v^{n-1}_{n-1}
\endbmatrix
$$
we have $\det \tilde V\neq 0$.

For $j=1,\ldots,n-1$ we find $t^j\in\Bbb R^{n-1}$ such that
$\tilde V t^j=e^j\in\Bbb R^{n-1}$.
Put $w^j:=\sum_{k=1}^{n-1}t^j_kv^k$, $j=1,\ldots,n-1$. We have
$w^j_l=\delta_{jl}$, $j,l=1,\ldots,n-1$.
Certainly, $w^j\in H$, $j=1,\ldots,n-1$.
It follows from the assumption of the lemma that the set
$\{w_n^1,\ldots,w_n^{n-1}\}$ is $\Bbb Z$-linearly independent
(that is if $\sum_{j=1}^{n-1}s_jw_n^j\in\Bbb Z$ for some
$s_j\in\Bbb Z$ then $s_j=0$). Then in
in view of multidimensional Kronecker Approximation Theorem
(see e.g. \cite{Hla-Sch-Tas}) the set
$$
\{(\alpha_n w_n^1-[\alpha_nw_n^1],\ldots,
\alpha_n w_n^{n-1}-[\alpha_nw_n^{n-1}]):
\alpha_n\in\Bbb Z\}
$$
is dense in $[0,1)^{n-1}$. But $\langle\alpha,w^j\rangle=\alpha_j+\alpha_n
w_n^j$; therefore,
$$
\{(\langle\alpha,w^1\rangle,\ldots,\langle\alpha,w^{n-1}\rangle):\alpha
\in\Bbb Z^n\}\quad\text{is dense in $\Bbb R^{n-1}$}.\tag{1}
$$

Put $T:=[t^1,\ldots,t^{n-1}]\in\Bbb R^{(n-1)\times (n-1)}$.
We have that $\det T\neq 0$. We have that
$$
[w^1,\ldots,w^{n-1}]=[v^1,\ldots,v^{n-1}] T.
$$
Consequently,
$$
(\langle\alpha,v^1\rangle,\ldots,\langle\alpha,v^{n-1}\rangle)=
(\langle\alpha,w^1\rangle,\ldots,\langle\alpha,w^{n-1}\rangle)T^{-1},
$$
which, in view of \thetag{1}, finishes the proof of the lemma.
\qed
\enddemo

\proclaim{Lemma 7 {\rm (cf. \cite{Zwo~3})}} Let $D$ be a bounded pseudoconvex domain
in $\Bbb C^n$. Fix $z^0\in\partial D$ satisfying the following condition
$$
\text{for any $j\in\{1,\ldots,n\}$ if $z_j^0=0$ then $D\cap V_j\neq\emptyset$}
$$
(this condition is satisfied if, for instance, $z^0\in\Bbb C_*^n$).

Then the condition \thetag{KC} is satisfied
at $z^0$ (for the subspace $\Cal E$).
\endproclaim
\demo{Proof} First note that for any $\alpha\in\Bbb Z^n$
such that $z^{\alpha}\in L_h^2(D)$ we have that $\alpha_j\geq 0$ 
if $z_j^0=0$. Therefore, it is sufficient to show that
$K_D(z)\to\infty$ as $z\to z^0$. Let $I:=\{j:z_j^0=0\}$. Without
loss of generality we may assume that $I=\{1,\ldots,s\}$.
We easily see that $s<n$. Then $D\subset\Bbb C^s\times\pi_I(D)$
(we identify $\pi_I(D)$ with a subset of $\Bbb C^{n-s}$, 
if $s=0$ then $\pi_I:=\id$). Note that the assumptions
of the criterion from 
\cite{Pfl~1} ('outer cone condition') are satisfied for the domain $\tilde D$
(and consequently also for $D$), where $\tilde D$ is a bounded
pseudoconvex Reinhardt domain in $\Bbb C^{n-s}$, $\pi_I(D)\subset\tilde D$,
$\pi_I(z^0)\in\partial \tilde D$ and $\partial\tilde D$ 
is $C^2$ near $\pi_I(z^0)$, which finishes the proof. The existence of such 
$\tilde D$ follows from the convexity of $\log \pi_I(D)$ and the fact that 
$\pi_I(z^0)\in\partial \pi_I(D)\cap\Bbb C_*^{n-s}$.
\qed
\enddemo
\proclaim{Lemma 8} Let $\beta,v\in\Bbb R^n$, $||v||=1$,
$\{x^{\nu}\}_{\nu=1}^{\infty}\subset\Bbb R^n$
be such that $||x^{\nu}||\to\infty$,  $\frac{x^{\nu}}{||x^{\nu}||}\to v$
as $\nu\to\infty$, and $\langle\beta,v\rangle<0$. Then
$$
\langle x^{\nu},\beta\rangle\to-\infty,\;\nu\to\infty.
$$
\endproclaim
\demo{Proof} Suppose that the Lemma does not hold. 
Then we may assume without loss of generality that
$\langle x^{\nu},\beta\rangle\geq M$ for some $M>-\infty$, $\nu=1,2,\ldots$.
Therefore,
$$
\langle \frac{x^{\nu}}{||x^{\nu}||},\beta\rangle\geq\frac{M}{||x^{\nu}||}.
$$
Passing with $\nu$ to infinity we get that $\langle v,\beta\rangle\geq 0$ --
contradiction.
\qed
\enddemo
Since $\frak C(D)\subset\Bbb R_-^n$, $\frak C(D)$ contains no straight
lines. The latter property is invariant with respect to linear isomorphisms
and is closely related to the hyperbolicity of a domain $D$ 
(see \cite{Zwo~1}); therefore, although it may 
be formulated a little more generally, we assume in Lemmas 9 and 10
that the cones contain no straight lines.

\proclaim{Lemma 9} Let $C$ be a
convex closed cone such that $C\cap\Bbb Q^n=\{0\}$
and let $C$ contain no straight lines. 
Let $v\in\wn_{\Span C}C$. Then 
for any $\delta>0$
there is $\beta\in\Bbb Z^n$ such that
$$
\gather
\langle\beta,v\rangle>0,\\
|\langle\beta,w\rangle|<\delta\;\text{ for any $w\in C$, $||w||=1$}.
\endgather
$$
\endproclaim
\demo{Proof} Denote by $\Cal U$ the largest vector subspace
of $\Span C$ among those spanned by vectors from $\Bbb Z^n$. 
Because of the assumptions of the lemma we have that $v\not\in\Cal U$.

Let $\{v^1,\ldots,v^r\}$ be a vector basis of $\Span C$ such that
$\{v^1,\ldots,v^s\}$ is a vector basis of $\Cal U$ and $v^r=v$. Certainly,
$s<r$. Since there is $M$ large enough such that
for any $w=\sum_{j=1}^rt_jv^j\in C$, $||w||=1$
we have that $|t_j|\leq M<\infty$, $j=1,\ldots,r$,
it is sufficient to find $\beta\in\Bbb Z^n$ 
such that $0\leq\langle\beta,v^j\rangle<\tilde\delta$, $j=1,\ldots,r$
and $\langle\beta,v^r\rangle>0$,
where $\tilde\delta:=\frac{\delta}{rM}>0$.

Let $A\in \Bbb Z^{n\times n}$ 
be a linear isomorphism of $\Bbb R^n$
such that $\Cal U=A(\Bbb R^s\times\{0\}^{n-s})$. Since
$\langle \gamma,Aw\rangle=\langle A^*\gamma,w\rangle$ 
for any $\gamma\in\Bbb Z^n$, $w\in\Bbb R^n$ we see that we may transportate
the problem to that with $\Cal U=\Bbb R^s\times\{0\}^{n-s}$
(possibly with other value of $\tilde\delta$).

Therefore, we assume that $\Cal U=\Bbb R^s\times\{0\}^{n-s}$.

Note that $\Span C\cap(\Bbb R^s\times\Bbb Q^{n-s})=
\Bbb R^s\times\{0\}^{n-s}$ and
the system of vectors $\{(v^j_{s+1},\ldots,v^j_n),
j=s+1,\ldots,r\}$ is linearly independent. 
Consequently, we get the existence of $\beta\in\Bbb Z^n$ such that
$\beta_j=0$, $j=1,\ldots,s$ and
$0<\langle\beta,v^j\rangle<\tilde\delta$, $j=s+1,\ldots,r$ 
(use Lemma 6 applied to $\pi_{1,\ldots,s}(\Span C)\subset
\{0\}^s\times\Bbb R^{n-s}$).
\qed
\enddemo

\proclaim{Lemma 10}
Let $C\subset\Bbb R^n$ be a convex closed cone 
such that $C\cap\Bbb Q^n=\{0\}$ and let $C$ contain no straight lines. Then 
for any $\delta>0$, $v\in C$, $v\neq 0$ 
there is $\beta\in\Bbb Z^n$ such that
$$
\gather
\langle\beta,v\rangle>0,\\
\langle\beta,w\rangle<\delta\; \text{ for any $w\in C$, $||w||=1$}.
\endgather
$$
\endproclaim
\demo{Proof} Denote by $H=H(C,v)$ a maximal vector subspace of $\Span C$ 
among those for which $v\in\wn_H(C\cap H)$ (one may verify that
$H$ is well-defined and unique). It follows from convexity of $C$ that
$$
\text{for any 
$w\in\wn_H(C\cap H)$ we have that $H(C,v)=H(C,w)$.}\tag{2}
$$
We shall need the following:
\proclaim{Claim 1} There are a sequence of vector subspaces $H=:H_k\subset
H_{k-1}\subset\ldots\subset H_1\subset H_0=\Span C$ such that
$\dim H_j=\dim H_{j+1}+1$, $j=0,\ldots,k-1$
and vectors $v^j\in H_j\setminus H_{j+1}$
orthogonal to $H_{j+1}$ such that
$$
C\cap H_j\subset\{u+tv^j:u\in H_{j+1},t\leq 0\}.\tag{3}
$$
\endproclaim
\demo{Proof of Claim 1} Note that to prove Claim 1 it is sufficient to show 
that having given a vector subspace $H_j$
of $\Span C$ such that $H\subsetneq H_j$ (or $j<k$) 
we can find a vector subspace $H_{j+1}$ with
$H\subset H_{j+1}\subset H_j$ and $\dim H_j=\dim H_{j+1}+1$,
and a vector $v^j$ such that \thetag{3} is satisfied.

There are two possibilities:

If $\Span(C\cap H_j)\neq H_j$ then we define $H_{j+1}$
as any vector hyperplane of $H_j$ containing $\Span(C\cap H_j)$.

If $\Span(C\cap H_j)=H_j$ then one may easily verify that 
$H\cap\wn_{H_j}(C\cap H_j)=\emptyset$ 
(it easily follows from definition of $H$ and \thetag{2}). Then
by the Hahn-Banach theorem there exists a supporting
hyperplane $H_{j+1}$ of $C\cap H_j$ in $H_j$ containing $H$, 
from which we easily get the desired $v^j$ and \thetag{3}.
\qed
\enddemo

It follows from Lemma 9 that there is $\tilde\beta\in\Bbb Z^n$ such that
$\langle\tilde\beta,v\rangle>0$ and $\langle\tilde\beta,w\rangle<\delta$ 
for any $||w||=1$, $w\in C\cap H$.
Therefore, applying induction and Claim 1 to finish the proof of Lemma 10
it is sufficient to prove the following:

\proclaim{Claim 2} Assume that there is 
$\tilde\beta$ as desired in Lemma 10
for $C\cap H_{j+1}$ and $v$ ($j<k$). 
Then there is $\beta\in\Bbb Z^n$ as desired in Lemma 10 for
$C\cap H_j$ and $v$ (the subspaces $H_j$ and $H_{j+1}$ are those appearing 
in Claim 1).
\endproclaim
\demo{Proof of Claim 2}
Put $M_1:=\sup\{\langle\tilde\beta,w\rangle:w\in H_{j+1},
||w||=1\}<\infty$. 

In view of the Dirichlet pigeon-hole theorem (see e.g. \cite{Har-Wri}) 
we have that for any positive integer $N$ there are
$\beta^N_l,q=q(N)\in\Bbb Z$, $q>0$ such that 
$\beta^N_l-qv^j_l=\epsilon(l,N)\in(-1/N,1/N)$, $l=1,\ldots,n$; moreover,
$q$ may be chosen so that it tends to infinity as $N$ tends to infinity.
Denote $\epsilon^N:=(\epsilon(1,N),\ldots,\epsilon(n,N))$. Then we have
$\beta^N=qv^j+\epsilon^N$.

We claim that for large $N$ $\beta:=\tilde\beta+\beta^N$ 
satisfies the desired property. 

First note that because $v^j$ is orthogonal to $v$ ($v\in H\subset H_{j+1}$)
we have
$$
\langle\tilde\beta+\beta^N,v\rangle=\langle\tilde\beta,v\rangle
+\langle\epsilon^N,v\rangle.
$$
Since the second summand in the formula above tends to $0$ 
as $N$ goes to infinity the last expression is positive for $N$ large enough.

Suppose that the second property of the lemma does not hold for infinetly
many $N$,
i.e. without loss of generality we may write that for any $N$ there are 
$$
t^N\leq 0,\;u^N\in H_{j+1},\; w^N=u^N+t^Nv^j
\in C,\; ||w^N||=1
$$ 
such that 
$\langle\tilde\beta+\beta^N,w^N\rangle\geq \delta$. 

Take $\tilde \delta<\delta$ such that $\langle\tilde\beta,u\rangle\leq\tilde 
\delta$, $u\in H_{j+1}\cap C$, $||u||=1$.

There is $M_2$ such that
for any $N$ we have $||u^N||\leq M_2$, $-M_2\leq t^N\leq 0$. 

Then we have ($v^j$ is orthogonal to $u^N$)
$$
\delta\leq\langle\tilde\beta+\beta^N,w^N\rangle=
\langle\tilde\beta,u^N\rangle+
t^N\langle\tilde\beta,v^j\rangle+\langle\epsilon^N,u^N\rangle+
t^N\langle\epsilon^N,v^j\rangle+t^Nq\langle v^j,v^j\rangle.\tag{4}
$$
Without loss of generality 
we may assume that $t^N\to t$. Moreover, we may assume that 
$u^N\to u\in H_{j+1}$ and, therefore, $w^N\to w\in C$, $||w||=1$. 

We claim that $t=0$. Suppose the contrary, so $t<0$. 
Then the first four summands
in the last expression of \thetag{4} are bounded from above
and the last expression tends to $-\infty$ -- contradiction.

Consequently, $w^N\to u$, $||u||=1$, $u\in C\cap H_{j+1}$. 
Note that in view of \thetag{4} 
(making use of the fact that $t^N\leq 0$, $q>0$), we have
$$
\delta\leq
\langle\tilde\beta,u^N\rangle+
t^N\langle\tilde\beta,v^j\rangle+\langle\epsilon^N,u^N\rangle+
t^N\langle\epsilon^N,v^j\rangle.
$$
Passing with $N$ to infinity we get 
$\delta\leq\langle\tilde\beta,u\rangle\leq \tilde\delta$ -- contradiction.
\qed
\qed
\enddemo
\enddemo

\subheading{3. Proof of Theorem 3}
The proof of implication
((i)$\implies$(ii)) is simple and is to find in \cite{Zwo~3}. 
For the sake of completeness
we give it here, too.
\demo{Proof of implication ((i)$\implies$(ii))}
Suppose that there is $v\in \frak C^{\prime}(D)\cap\Bbb Q^n$. Certainly, 
$v\neq 0$.
We assume that $a\in\log D$ 
from the definition of $\frak C(D)$ is
equal to $(0,\ldots,0)$.
Without loss of generality
we may assume that $v\in\Bbb Z_-^n$ and $v_1,\ldots,v_n$ are relatively
prime.

It is sufficient to show that the Bergman length $L_{\beta_D}$ of the
curve $(t^{-v_1},\ldots,t^{-v_n})$, $0<t<1$ is finite.

Denote $\phi(\lambda):=(\lambda^{-v_1},\ldots,\lambda^{-v_n})$,
$\lambda\in E_*$. Certainly, $\phi\in\Cal O(E_*,D)$.
Put $u(\lambda):=K_D(\phi(\lambda))$. Then we have (use Lemma 5)
$$
u(\lambda)=\sum_{\alpha\in\Bbb Z^n:\langle\alpha+\bold 1,v\rangle<0}
a_{\alpha}|\lambda|^{-2\langle\alpha,v\rangle}=\sum_{j=j_0}^{\infty}b_j
|\lambda|^{2j},
$$
where $b_{j_0}\neq 0$ (note that $j_0>\langle\bold 1,v\rangle$ and
it is possible that many of $b_j$'s in the formula above vanish).

Note that 
$$
\beta_D^2(\phi(\lambda);\phi^{\prime}(\lambda))=\frac{\partial^2
\log u(\lambda)}{\partial\lambda\partial\bar\lambda}=\frac{\partial^2}
{\partial\lambda\partial\bar\lambda}\left(\log\sum_{j=j_0}^{\infty}b_j
|\lambda|^{2j-2j_0}\right).
$$ 
The last expression tends to some constant $C\in\Bbb R$, 
which finishes the proof.
\qed
\enddemo
\demo{Proof of implication ((ii)$\implies$(i))}
We prove that the condition \thetag{KC}
is satisfied in all $z^0\in\partial D$ for $\Cal E$. 
Recall that we already know that \thetag{KC} is satisfied for all points 
$z^0\in\partial D$ such that
$$
\text{ for any $j$ (if $z_j^0=0$ then $V_j\cap D
\neq\emptyset$)}\tag{5}
$$
(see Lemma 7); in particular, 
for all points from $\partial D\cap\Bbb C_*^n$. 

Additionally, note that to prove the property \thetag{KC} (at $z^0$)
it is sufficient to consider sequences 
$\{z^{\nu}\}_{\nu=1}^{\infty}\subset D\cap\Bbb C_*^n$.

Without loss of generality $(1,\ldots,1)\in D$.
Take some $v\in\frak C(D)\cap\Bbb Q^n$, $v\neq 0$.
Since $v\in\tilde\frak C(D)$ (assumption of the theorem), 
we get from definition of $\tilde\frak C(D)$ that 
$$
\lim_{t\to\infty}(\exp(tv_1),\ldots,\exp(tv_n))=:w\in D.
$$
Note that $w_j=0$ if $v_j<0$ and $w_j=1$ if $v_j=0$. In particular,
$$
\text{for any $v\in\frak C(D)\cap\Bbb Q^n$ 
(if $v_j<0$ then $D\cap V_j\neq\emptyset$).}\tag{6}
$$
Without loss of generality we may assume that $D\cap V_j\neq\emptyset$,
$j=1,\ldots,k$, $D\cap V_j=\emptyset$, $j=k+1,\ldots,n$. Because we are
interested in these $z^0$ for which \thetag{5} 
is not satisfied we may assume that $k<n$.

In view of our assumptions 
(and properties of pseudoconvex Reinhardt domains) we know that 
$$
\Bbb R_-^k\times\{0\}^{n-k}\subset\frak C(D).\tag{7}
$$
Note that 
$$
\text{for any $v\in \frak C(D)\setminus(\Bbb R^k\times\{0\}^{n-k})$
we have that $v\not\in\Bbb R^k\times\Bbb Q^{n-k}$.}\tag{8}
$$
Actually, suppose that there exists
$v\in \frak C(D)\cap(\Bbb R^k_-\times\Bbb Q^{n-k}_-)$, 
$v_j<0$ for some $j>k$. Then adding some vector from 
$\Bbb R_-^k\times\{0\}^{n-k}$ we get a vector 
(we denote it with the same letter) from $\frak C(D)\cap\Bbb Q^n$ (with
$v_j<0$ for some $j>k$).
Then in view of \thetag{6} $D\cap V_j\neq\emptyset$ -- contradiction. 

Denote $\tilde x:=\pi(x)$, where $\pi(x):=(0,\ldots,0,x_{k+1},\ldots,x_n)$,
$x\in\Bbb R^n$. 
Note that $\pi(\frak C(D))$ is a closed convex cone 
in $\{0\}^k\times\Bbb R^{n-k}_-$ and in view of \thetag{8}
$\pi(\frak C(D))\cap(\{0\}^k\times\Bbb Q^{n-k})=\{0\}$.

Consider a point $z^0\in\partial D$, not satisfying \thetag{5}, 
in particular,
$$
\text{there is $j>k$ such that $z^0_j=0$,}\tag{9}
$$ 
and a sequence
$z^{\nu}\to z^0$, $z^{\nu}\in D\cap\Bbb C_*^n$. 
Put $x^{\nu}:=\log|z^{\nu}|$.  Without loss
of generality we may assume that $\frac{x^{\nu}}{||x^{\nu}||}\to v\in 
\frak C(D)$. Certainly, $||x^{\nu}||\to\infty$.

Fix $\alpha\in\Bbb Z^n$ such that $z^{\alpha}\in L_h^2(D)$. 
Define $\delta:=\inf\{-\langle\alpha+\bold 1,w\rangle:w\in\frak C(D),
||w||=1\}>0$ (use Lemma 5).

Below we consider two cases:

Case (I). $v_j<0$ for some $j>k$. 

We claim that it is sufficient to find $\beta\in\Bbb Z^n$ such that
$$
\langle\beta,w\rangle<\delta \text{ for any $w\in\frak C(D)$, $||w||=1$
and $\langle\beta,v\rangle>0$}.\tag{10}
$$ 
In fact, then $z^{\alpha+\beta}\in L_h^2(D)$ (use Lemma 5) and
$$
\frac{|(z^{\nu})^{\alpha}|}{\sqrt{K_D(z^{\nu})}}\leq
||z^{\alpha+\beta}||_{L^2(D)}
\frac{|(z^{\nu})^{\alpha}|}
{|(z^{\nu})^{\alpha+\beta}|}
=||z^{\alpha+\beta}||_{L^2(D)}
|(z^{\nu})^{-\beta}|.
$$
The last expression tends to zero (use Lemma 8).

Therefore, we prove the existence of $\beta\in\Bbb Z^n$
such that \thetag{10} is satisfied.

Use Lemma 10 (applied to $\pi(\frak C(D))$ and $\pi(v)$)
to get the existence of $\beta\in\{0\}^k\times\Bbb Z^{n-k}$
such that $\langle\beta,v\rangle=\langle\beta,\pi(v)\rangle>0$
and $\langle\beta,w\rangle=||\pi(w)||\langle\beta,\frac{\pi(w)}{||\pi(w)||}
\rangle<\delta$
for any $w\in\frak C(D)$, $||w||=1$ with $\pi(w)\neq 0$. 
Since $\langle\beta,w\rangle=0$
if $\pi(w)=0$, this finishes the proof.

Case (II). $v_{k+1}=\ldots=v_n=0$.

Put $\tilde x:=\pi(x)$.
Without loss of generality 
$\frac{\tilde x^{\nu}}{||\tilde x^{\nu}||}\to\tilde w\in\{0\}^k\times
\Bbb R^{n-k}$. In view of \thetag{9} we have $||\tilde x^{\nu}||\to\infty$.

Note that it is sufficient to find $\beta\in\{0\}^k\times\Bbb Z^{n-k}$
such that $\langle\beta,\tilde w\rangle>0$ and $\langle\beta,w\rangle=
||\pi(w)||\langle\beta,
\frac{\pi(w)}{||\pi(w)||}\rangle<\delta$, where $\delta$ is as earlier,
$w\in\frak C(D)$, $w\neq 0$.
Then similarly as earlier we have that $z^{\alpha+\beta}\in L_h^2(D)$
and
$$
\frac{|(z^{\nu})^{\alpha}|}{\sqrt{K_D(z^{\nu})}}\leq
||z^{\alpha+\beta}||_{L^2(D)}
|(z^{\nu})^{-\beta}|=
||z^{\alpha+\beta}||_{L^2(D)}
|(z^{\nu}_{k+1},\ldots,z^{\nu}_n)^{-(\beta_{k+1},\ldots,\beta_n)}|.
$$
The last expression tends to zero (remember about convergence 
$||\tilde x^{\nu}||\to\infty$ and then use Lemma 8).

If $\tilde w\in\pi(\frak C(D))$ then the existence of such $\beta$ 
follows from Lemma 10 applied to $\pi(\frak C(D))$ and $\tilde w$.

If $\tilde w\not\in\pi(\frak C(D))$ 
then define $\tilde C$ to be the smallest convex
cone containing $\frak C(D)$ and $-\tilde w$.
It is easy to see that 
$\tilde C\neq\{0\}^k\times\Bbb R^{n-k}$ (e.g. $\tilde w\not\in\tilde C$). 
Consequently, the set
$$
\{\beta\in\{0\}^k\times\Bbb R^{n-k}:\langle\beta,u\rangle<0,
u\in\tilde C\setminus\{0\}\}
$$
is a non-empty convex open cone (in $\{0\}^k\times\Bbb R^{n-k}$). 
Therefore, it contains 
$\beta\in\{0\}^k\times\Bbb Z^{n-k}$. In particular,
$\langle\beta,-\tilde w\rangle<0$, $\langle\beta,\frac{\pi(w)}{||\pi(w)||}
\rangle<0<\delta$
for any $w\in \frak C(D)$, $\pi(w)\neq 0$,
which finishes the proof.
\qed
\enddemo

\subheading{4. Concluding remarks} In view of the results from \cite{Pfl~2},
\cite{Fu} and \cite{Zwo~2} we may precisely give the relations between
hyperconvexity, Carath\'eodory completeness, Kobayashi completeness
and Bergman completeness in the class of bounded pseudoconvex Reinhardt domains
and we may express it with the help of the set
$\frak C^{\prime}(D)$. Although there are many not biholomorphic bounded
pseudoconvex Reinhardt domains with the same $\frak C(D)$ 
and $\frak C^{\prime}(D)$ this is the set $\frak C^{\prime}(D)$, 
which desribes entirely the completeness of the domain $D$. 

Summarizing, the following properies are satisfied:

--- all bounded Reinhardt pseudoconvex domains are Kobayashi complete;

--- hyperconvexity is equivalent to Carath\'eodory completeness
and the last is equivalent to the equality $\frak C^{\prime}(D)=\emptyset$;

--- Bergman completeness is equivalent to the equality
$\frak C^{\prime}(D)\cap\Bbb Q^n=\emptyset$.

In view of the above remarks we may easily produce a great variety
of bounded Reinhardt domains, which are Bergman complete but not hyperconvex.

Note that the proper choice of the subspace $\Cal E$ and the Kobayashi
criterion reduce the problem of the proof of Bergman completeness of the
considered domains to a problem from diophantine approximation on 
cones $\frak C(D)$.

\subheading{Acknowledgments} The paper was written after many stimulating 
discussions with professor Peter Pflug. The author would like to thank him.

\enddocument